%% file: main.tex
\documentclass[11pt, number, times, preprint]{article}
\usepackage[margin=1.7in]{geometry}

\usepackage{amsmath,amsfonts,amssymb,amsthm,epsfig,epstopdf,titling,url,array, pdflscape, fancyhdr, lastpage, tikz, hyperref}
\usepackage{blindtext}
\usepackage{dsfont}
\usepackage{caption}
\usepackage{subcaption}

\title{Skew Axial Algebras of Monster Type}
\author{Michael Turner\textsuperscript{}\thanks{\textsuperscript{}
School of Mathematics, University of Birmingham, Edgbaston, Birmingham, B15 2TT, UK, Email:
mxt187@student.bham.ac.uk
}}
\date{\today}

\newtheoremstyle{sltheorem}
{}                
{}                
{\slshape}        
{}                
{\bfseries}       
{.}               
{ }               
{}                

\theoremstyle{sltheorem}
\newtheorem{thm}{Theorem}[section]
\newtheorem{lem}[thm]{Lemma}
\newtheorem{prop}[thm]{Proposition}
\newtheorem{cor}[thm]{Corollary}
\newtheorem*{thm*}{Theorem}
\newtheorem*{cor*}{Corollary}

\theoremstyle{definition}
\newtheorem{defn}[thm]{Definition}

\newtheorem{exmp}[thm]{Example}

\newtheorem{no}[thm]{Notation}

\theoremstyle{remark}
\newtheorem*{rem}{Remark}
\newtheorem*{note}{Note}

\newcommand{\id}{\mathds{1}}
\newcommand{\F}{\mathbb{F}}
\newcommand{\mcal}[1]{\mathcal{#1}}

\newcommand{\GenG}[1]{\langle #1\rangle}
\newcommand{\GenA}[1]{\langle\!\langle #1\rangle\!\rangle}
\newcommand{\fustar}[1]{(\mcal{#1},\star)}

\newcommand{\Z}{\mathbb{Z}}
\newcommand{\N}{\mathbb{N}}

\newcommand{\spec}[1]{\text{Spec}(#1)}

\newcommand{\Span}[1]{\text{span}(#1)}

\newcommand{\al}{\alpha}
\newcommand{\bt}{\beta}
\newcommand{\sg}{\sigma}
\newcommand{\lm}{\lambda}
\newcommand{\dt}{\delta}
\newcommand{\ep}{\epsilon}
\newcommand{\lmf}{\lambda^f}
\newcommand{\tu}[1]{\tau_{#1}}
\newcommand{\gm}{\gamma}
\newcommand{\te}{\theta}
\newcommand{\zt}{\zeta}
\newcommand{\kp}{\kappa}

\newcommand{\jor}[1]{\mcal{J}(#1)}
\newcommand{\mon}[1]{\mcal{M}(#1)}

\newcommand{\TB}{2\text{B}}

\newcommand{\TC}{3\text{C}}

\begin{document}
\maketitle
\input{Abstract}
\input{Introduction}
\input{Prerequisites}
\input{Examples}
\input{Construction}
\input{SkewRelations}

\input{TheoremProof}

\input{CorollaryProof}
\input{Consequences}
\input{Appendix}
\bibliographystyle{abbrv}
\bibliography{references}
\end{document}

%% file: Abstract.tex
\begin{abstract}
Skew axets were first defined by M\textsuperscript{c}Inroy and Shpectorov where they used the term of axets to classify shapes of an algebra. When they first submitted their paper, it was not known if skew axial algebras exist and now we will present such examples with axet $X'(1+2)$. Looking at $2$-generated primitive axial algebras of Monster type, we will be able to state and prove the classification of such algebras with axet $X'(1+2)$. We will conclude by looking at larger skew axets and give a suggestion on how one could extend the classification.
\end{abstract}

%% file: Introduction.tex
\section{Introduction}
Hall, Shpectorov and Rehren first defined axial algebras in \cite{hall2015primitive} and \cite{hall2015universal}. These are non-associative, commutative algebras which are strongly related to the $3$-transposition groups as well as some of the sporadic groups. They were inspired by Majorana algebras which were defined in \cite{ivanov2009monster} using work around the properties of the Griess algebra. In \cite{ivanov2010majorana}, the $2$-generated Majorana algebras were classified. The properties of Majorana algebras can be generalised even further to the other fusion laws, for example $\mon{\al,\bt}$, to produce what are called axial algebras.

To motivate why skew axial algebras are a topic of interest, we look at recent work of M\textsuperscript{c}Inroy and Shpectorov in \cite{mcinroy2021forbidden}. Their paper clarified a problem with shapes of an algebra. Shapes were defined from algebras however there are given shapes where no algebra exists, resulting in a circular argument. Their goal was to make the definition of shape independent of algebras to avoid this oxymoron. While tackling this challenge, they defined axets using group theoretic properties rather than through axial algebras. They noticed in a certain case that these axets are either regular or skew. As the known axial algebras are regular, they posed a question of if any axial algebra is skew? In this paper, we will answer such question.

As there are an infinite number of skew axets, we will be focusing on the smallest one, $X'(1+2)$. This axet does have axial algebra examples, which we will discuss later on. We will start by briefly introducing axial algebras and axets to the reader. This will not be a deep discussion and it will keep to the essentials. We will then give three examples of skew axial algebras. The first is $\TC(\al)$ when we assign a Monster fusion law to this algebra rather than the commonly assigned Jordan fusion law. For $\al=-1$, we change our approach to get a skew axial algebra with fusion law $\mon{-1,2}$. The second will be $Q_2(\frac{1}{3})$, an algebra which is related to Matsuo algebras. Again, this algebra is assigned a different fusion law $\mon{\frac{1}{3},\frac{2}{3}}$ rather than $\mon{\frac{2}{3},\frac{1}{3}}$. A third example, $Q_2(\frac{1}{3})^\times \oplus \GenG{\id}$ is also constructed in characteristic $5$. This is due to $Q_2(\frac{1}{3})$ being not simple in characteristic $5$ and it has a $3$-dimensional quotient, $Q_2(\frac{1}{3})^\times$. Adding a universal identity element, we get a new $2$-generated axial algebra which is also skew and with a Monster fusion law. This gives an answer to Problem 6.14 in \cite{mcinroy2022axial} however a complete classification is still in progress.

Convincing the reader that these skew algebras exist, the rest of the paper will be focused on classifying these algebras for skew axet $X'(1+2)$. This will be done by applying work from \cite{rehren2017generalised} and \cite{franchi20211} to produce a complete multiplication table. In the appendix, we will produce relations between constants, using GAP \cite{gap2022}, to make our lives easier. With this information, the proof is split into two parts:
\begin{itemize}
\item when two of the axes are orthogonal, and
\item when they are not orthogonal. 
\end{itemize}
When axes are orthogonal, this gives rise to double axes, which has been studied in multiple papers, for example \cite{joshi2020axial}. Looking at those parts separately, our work proves the following:
\begin{thm}\label{Theorem}
Let $\F$ be a field of characteristic not equal to $2$. Suppose $(A,X)$ is a $2$-generated primitive $\mon{\al,\bt}$-axial algebra over $\F$ and $\GenG{X}$ is isomorphic to $X'(1+2)$. Then either
\begin{enumerate}
\item[$1.$] $A \cong \TC(\al,1-\al)$ for $\al\neq \frac{1}{2}$ and $\al+\bt=1$, or
\item[$2.$] $(\al,\bt)=(\frac{1}{3},\frac{2}{3})$ and either
\begin{enumerate}
    \item[$i)$] $A\cong Q_2(\frac{1}{3},\frac{2}{3}$) if $\F$ has characteristic not equal to $5$, or
    \item[$ii)$] $A\cong Q_2(\frac{1}{3})^\times\oplus \GenG{\id}$ if $\F$ has characteristic equal to $5$.
\end{enumerate}
\end{enumerate}
\end{thm}
\begin{rem}
Adopting the notation in \cite{franchi20212}, we will denote the algebra $\TC(\al)$, $\al\neq -1$, with Monster fusion law $\mon{\al,1-\al}$ as $\TC(\al,1-\al)$ to avoid any confusion about which fusion law we have on this algebra. We denote $\TC(-1,2)$ to be $\TC(2)$ with $\mon{-1,2}$ fusion law.  Further, $Q_2(\frac{1}{3},\frac{2}{3})$ is $Q_2(\frac{1}{3})$ with a different fusion law to its definition in \cite{galt2021double}. The final algebra is defined in Table \ref{mult Qx}.
\end{rem}
\begin{note}
In all cases, $\al+\bt=1$. This is due to one of the axes in the axet being twisted by the identity element.
\end{note}
We will then prove the following corollary. 
\begin{cor}\label{Corollary}
Let $\F$ be a field of characteristic not equal to $2$. Suppose that $A=\GenA{p,q}$ is a primitive axial algebra over $\F$ such that $p$ is a $\mon{\al,\bt}$-axis and $q$ is a $\mcal{J}(\al)$-axis. Then $A$ is isomorphic to one of the following:
\begin{enumerate}
    \item[$1.$] an axial algebra of $\jor{\al}$-type, or
    \item[$2.$] an axial algebra listed in the Theorem $\ref{Theorem}$.
\end{enumerate}
\end{cor}

The final section will start the discussion on larger skew axets, $X'(k+2k)$, and how our result could help with such task. When $k$ is odd and choosing the right axes, there will be a subaxet of $X'(1+2)$ inside the larger skew axet. We will conclude by proposing a plan to tackle the odd $k$ case, which we hope to present in a secondary paper with work on even $k$.

The author would like to thank the anonymous referee for their useful comments and corrections.

%% file: Prerequisites.tex
\section{Prerequisites}
\subsection{Axial Algebras}
Axial algebras have been defined in multiple papers so this section will be brief. For an in-depth discussion, we recommend to look at \cite{hall2015universal}, \cite{khasraw2020structure} and \cite{mcinroy2022axial}. A starting point is the definition of a fusion law in \cite{de2020decomposition}. 
\begin{defn}
Let $\mcal{F}$ be a set and let $\star: \mcal{F} \times \mcal{F} \rightarrow 2^{\mcal{F}}$ be a binary operator, where $2^{\mcal{F}}$ denotes the power set of $\mcal{F}$. We call the pair $\fustar{F}$ a \emph{fusion law}. We say it is  \emph{symmetric} if for all $x,y\in \mcal{F}$, we have $x\star y=y\star x$.  
\end{defn}
\begin{note}
When a fusion law $\fustar{F}$ is clear, we will refer to it by $\mcal{F}$. 
\end{note}
For axial algebras, as they are commutative, fusion laws are considered to be symmetric. Given a fusion law, one is able to grade it using a group. We will denote the group operation by juxtaposition. 
\begin{defn} Let $\fustar{F}$ be a fusion law and $T$ be a group. We say $\mcal{F}$ is $T$-graded if there exists a map $\xi: \mcal{F}\rightarrow T$
such that for all $\lm, \mu\in \mcal{F}$,
$$\xi(\nu)=\xi(\lm)\xi(\mu)$$
for all $\nu\in \lm\star \mu$.
\end{defn}
\begin{table}[b]
\begin{minipage}{.5\linewidth}
\centering
\begin{tabular}{|c|c|c|c|} 
 \hline
  $\star$& 0 & 1 & $\eta$\\ [0.4ex] 
 \hline
 0 & 0 &  & $\eta$\\ 
 \hline
 1 &  & 1 & $\eta$\\
 \hline
 $\eta$ & $\eta$ & $\eta$& 1, 0\\
 \hline
\end{tabular}
\caption{Fusion law of $\jor{\eta}$}\label{J Table}
\end{minipage}%
\begin{minipage}{.5\linewidth}
\centering
\begin{tabular}{|c|c|c|c|c|} 
 \hline
  $\star$& 0 & 1 & $\alpha$ & $\beta$\\ [0.4ex] 
 \hline
 0 & 0 &  & $\alpha$ & $\beta$\\ 
 \hline
 1 &  & 1 & $\alpha$ & $\beta$\\
 \hline
 $\alpha$ & $\alpha$ & $\alpha$& 1, 0 & $\beta$\\
 \hline
 $\beta$ & $\beta$ & $\beta$ & $\beta$ & 1, 0, $\alpha$\\
 \hline
 \end{tabular}
 \caption{Fusion law of $\mon{\al,\bt}$}
 \label{M Table}
 \end{minipage}
\end{table}
The fusion laws that are used in this paper are $\jor{\eta}$ and $\mon{\al,\bt}$, which are defined in Tables \ref{J Table} and \ref{M Table} respectively. Let $C_2=\{e,s\}$ be the cyclic group of order two, where $e$ is the identity element. Then $\jor{\eta}$ and $\mon{\al,\bt}$ are both $C_2$-graded where $0$, $1$, $\al$ are mapped to $e$ and $\eta$, $\bt$ are mapped to $s$.

Let $A$ be a commutative algebra over $\F$. For $a\in A$, the adjoint map of $a$ is $\text{ad}_a: A \rightarrow A$ with $\text{ad}_a(x)=ax$ for all $x\in A$. Let $\text{Spec}(a)$ be the set of eigenvalues of $\text{ad}_a$ and for $\lambda \in \F$, let $A_{\lambda}(a)$ be the $\lambda$-eigenspace of $\text{ad}_a$. Note that $A_\lambda(a)$ is non-trivial if and only if $\lambda \in \spec{a}$. For $S\subseteq \F$, we write $A_S(a) := \oplus_{\lm \in S} A_\lm(a)$.
We will use the definition of axes and axial algebras in \cite{khasraw2020structure}. 
\begin{defn}
Let $\fustar{F}$ be a fusion law, $A$ be a commutative algebra over $\F$, and $a\in A$. We have that $a$ is a (primitive) $\mcal{F}$-\emph{axis} if it satisfies the following:
\begin{enumerate}
    \item[A$1.$] $a$ is an idempotent; $a^2=a$,
    \item[A$2.$] $a$ is semisimple and $\spec{a}\subseteq \mcal{F}$; that is, $A=A_\mcal{F}(a)$,
    \item[A$3$.] For all $\lm,\mu \in \mcal{F}$, $A_\lm(a)A_\mu(a) \subseteq A_{\lm \star \mu}(a)$, and 
    \item[A$4$.] $A_1(a)=\GenG{a}$.
\end{enumerate}
We call an $\mcal{F}$-axis non-primitive if it satisfies A$1$-A$3$ and not A$4$. 
\end{defn}
\begin{defn}
Let $\fustar{F}$ be a fusion law, $A$ be a commutative algebra over $\F$, and $X$ be a set of (primitive) $\mcal{F}$-axes in $A$. We call $(A,X)$ a (primitive) $\mcal{F}$-\emph{axial algebra} if $X$ generates $A$.
\end{defn}
\begin{no}
When $X$ and $\mcal{F}$ are clear, we will refer to the axial algebra by $A$. As this paper is concerned with only $2$-generated algebras, we will assume $|X|=2$. We will assume all axial algebras and axes to be primitive. If that is not the case, we will explicitly say when an axis or an axial algebra is not.  
\end{no}
\begin{defn}\label{Proj}
    Let $(A,X)$ be a $\mcal{F}$-axial algebra and let $a$ be a $\mcal{F}$-axis. For any $v\in A$, $v=\lm_a(v) a + \sum_{\gm \in \mcal{F}\backslash\{1\}} v_\gm$ where $v_\gm \in A_\gm(a)$ due to A$2$. Let
    \begin{eqnarray*}
        \lm_a : A &\rightarrow& \F\\
 v &\mapsto& \lm_a(v).      
    \end{eqnarray*}
    We call this map the \emph{projection map} of $a$. This is a well-defined $\F$-linear map from Proposition 2.4 in \cite{franchi20211}.
\end{defn}
\begin{defn}
Let $\F$ have characteristic not equal to $2$. Suppose $\fustar{F}$ is a $C_2$-graded fusion law with morphism $\xi$, and $(A,X)$ is a $\mcal{F}$-axial algebra. Suppose $a$ is a $\mcal{F}$-axis in $A$ and let $\tu{a}:A\rightarrow A$. For $v\in A_\lm(a)$, we define $\tu{a}$ to be the following 
\begin{eqnarray*}
        \tau_a(v)=
        \begin{cases}
            v & \text{if } \xi(\lm)=e\\
            -v & \text{if } \xi(\lm)=s
        \end{cases}
    \end{eqnarray*}
By A$2$, this definition of $\tu{a}$ extends to the entire algebra and it can be shown to an automorphism of $A$. We call $\tu{a}$ the \emph{Miyamoto involution} of $a$.
\end{defn}
\begin{exmp}
For $(A,X)$ a $\jor{\eta}$-axial algebra and $\jor{\eta}$-axis $a$, $\tu{a}$ acts as the identity map on $A_{\{0,1\}}(a)$ and negative identity map on $A_\eta(a)$.

Similarly for $(A,X)$ a $\mon{\al,\bt}$-axial algebra and $\mon{\al,\bt}$-axis $a$, $\tu{a}$ acts as the identity map on $A_{\{0,1,\al\}}(a)$ and negative identity map on $A_\bt(a)$.
\end{exmp}

\begin{defn}
We call a fusion law, $\fustar{F}$, \emph{Seress} if:
\begin{enumerate}
    \item $0 \in \mcal{F}$, and
    \item $\lambda \star 0 = \{\lambda\}$ for $\lambda \neq 1$, and $1\star 0 = \emptyset$.
\end{enumerate}
\end{defn}
\begin{exmp}
Both $\jor{\eta}$ and $\mon{\al,\bt}$ are Seress.
\end{exmp}
The following result was proven in \cite{hall2015universal}. Seress originally proved it for fusion law $\mon{\frac{1}{4},\frac{1}{32}}$ but has now been generalised to all fusion laws that are Seress. 
\begin{lem}[Seress Lemma]\label{Seress}
Let $\fustar{F}$ be Seress and $(A,X)$ be an $\mcal{F}$-axial algebra. Then every $\mcal{F}$-axis $a\in A$ associates with $A_{\{0,1\}}(a)$. That is, for $x\in A$ and $u\in A_{\{0,1\}}(a)$, we have
\[ a(xu)=(ax)u.\]
\end{lem}
This result will be useful when finding the GAP equations in the appendix.
\subsection{Axets}
The idea of axets can be traced back to \cite{mcinroy2020expansion}, where M\textsuperscript{c}Inroy and Shpectorov looked at shapes of an algebra for their algorithm. Later in \cite{mcinroy2021forbidden}, they defined axets and stated results around them. We recommend that the reader looks at \cite[Section 3]{mcinroy2021forbidden} and \cite[Section 6.3]{mcinroy2022axial} for a greater introduction.
\begin{defn}
Let $S$, $G$ be groups and $G$ acts on a set $X$. Suppose there is a map $\tau:X\times S \rightarrow G$, denoted $\tau(x,s)=\tu{x}(s)$. Then $(G,X,\tau)$ is called an $S$-axet if for all $x\in X$, $s,s'\in S$ and $g\in G$, the following properties hold:
\begin{enumerate}
    \item $\tu{x}(s)\in G_x$;
    \item $\tu{x}(ss')=\tu{x}(s)\tu{x}(s')$; and 
    \item $\tu{xg}(s)=\tu{x}(s)^g$.
\end{enumerate}
If $S$, $G$ and $\tau$ are clear, we denote the axet by $X$.
For each $x\in X$, let $T_x:=\text{Im}(\tu{x})$, which is called the axial subgroup corresponding to $x$.
\end{defn}
\begin{rem}
The first two points makes $\tau_x$ a group homomorphism between $S$ and $G_x$ for all $x\in X$. One can show that for all $x\in X$, $T_x\subseteq Z(G_x)$ and $[S,S]\subseteq \text{ker}(\tu{x})$. Hence we can assume $S$ to be abelian.   
\end{rem}

\begin{defn}
Let $(G,X,\tau)$ be an $S$-axet and $Y\subseteq X$. Then $Y$ is \emph{closed} if $Y$ is invariant under $T_y$ for all $y\in Y$. For $Z\subseteq X$, the \emph{closure} of $Z$ in $X$ is the smallest closed subset containing $Z$, denoted by $\GenG{Z}$. We say that $\GenG{Z}$ is $k$-\emph{generated} if there exist $x_1,...,x_k\in X$ such that $\GenG{Z}=\GenG{x_1,...,x_k}$.
\end{defn}
\begin{exmp}
Let $(A,X)$ be an $\mcal{F}$-axial algebra where $\mcal{F}$ is $C_2$-graded. Take $\tilde{X}$ to be the set of all $\mcal{F}$-axes in $A$ and $G=\text{Aut}(A)$. For $x\in \tilde{X}$, let $\tu{e}(x)$ to be the identity map and $\tu{s}(x)$ to be the Miyamoto involution of $x$. Then $(G,\tilde{X},\tau)$ is a $C_2$-axet.

Let $Y=\GenG{X}$, $H=G_Y$ and $\tau'$ be the restriction of $\tau$ to the domain $Y\times C_2$. Then $(H,Y,\tau')$ is $C_2$-axet. Suppose that $|X|=2$ whence $(H,Y,\tau')$ is a $2$-generated axet and this paper will solely be looking at these.
\end{exmp}
We will now discuss some $2$-generated $C_2$-axets. 
\begin{exmp}
The following are in Examples 3.22 and 3.25 in \cite{mcinroy2021forbidden}.
\begin{enumerate}
    \item For $n\geq 2$, denote $X(n)$ to be the vertices of a regular $n$-gon. Let $X=X(n)$ and $a$, $b$ be adjacent vertices in $X$. We have that $G=D_{2n}$ and for $a\in X$, $\tu{a}(e)$ is the identity map and $\tu{a}(s)$ is the reflection of $X$ at vertex $a$. 
    \item Let $X=X(\infty)\cong\Z$ and $a$, $b$ be adjacent points on $X$. Then $G=D_\infty$ and for $a\in X$, $\tu{a}(e)$ is the identity map and $\tu{a}(s)$ is the reflection of $X$ at point $a$.
    \item Let $n=4k$ and we look at the $C_2$-axet $(D_{8k}, X(4k), \tau)$. We now identify opposite vertices of one bipartite half of the polygon denoting this by $X':=X'(k+2k)$. Let $G=D_{4k}$ and $\tau'$ to be the restriction of the domain to $X'\times C_2$. This is also an $C_2$-axet and we call $X'$ \emph{skew}. 
\end{enumerate}
\end{exmp}
We will now state Theorem 3.27 in \cite{mcinroy2021forbidden}.
\begin{prop}\label{axet thm}
Let $X=\GenG{a,b}$ be a $2$-generated $C_2$-axet with $n$ axes, where $n\geq 2$ and $n=\infty$. Then $X$ is isomorphic to 
\begin{enumerate}
    \item[$1.$] $X(n)$, or
    \item[$2.$] $X'(k+2k)$, where $k\in \N$ and $n=3k$.
\end{enumerate}
\end{prop}
In Figures \ref{Skew-1} and \ref{Skew-2}, we show $X'(1+2)$ and $X'(2+4)$, where the dotted arrow represents equivalence between two vertices. From now on, we will be working with $X'(1+2)$ however we will talk about more general skew axets in the final section. 

\begin{figure}[ht]
\centering
\begin{minipage}{.5\textwidth}
\centering
\tikzset{every picture/.style={line width=0.75pt}} 

\begin{tikzpicture}[x=0.75pt,y=0.75pt,yscale=-1,xscale=1]

\draw   (307.01,75.59) -- (386.35,154.93) -- (307.01,234.28) -- (227.66,154.93) -- cycle ;
\draw [color={rgb, 255:red, 255; green, 4; blue, 8 }  ,draw opacity=1 ] [dash pattern={on 0.84pt off 2.51pt}]  (307.01,78.59) -- (307.01,231.28) ;
\draw [shift={(307.01,234.28)}, rotate = 270] [fill={rgb, 255:red, 255; green, 4; blue, 8 }  ,fill opacity=1 ][line width=0.08]  [draw opacity=0] (8.93,-4.29) -- (0,0) -- (8.93,4.29) -- cycle    ;
\draw [shift={(307.01,75.59)}, rotate = 90] [fill={rgb, 255:red, 255; green, 4; blue, 8 }  ,fill opacity=1 ][line width=0.08]  [draw opacity=0] (8.93,-4.29) -- (0,0) -- (8.93,4.29) -- cycle    ;

\draw (301,60) node [anchor=north west][inner sep=0.75pt]   [align=left] {$\displaystyle a$};
\draw (212,147) node [anchor=north west][inner sep=0.75pt]   [align=left] {$\displaystyle b$};
\draw (301,240) node [anchor=north west][inner sep=0.75pt]   [align=left] {$\displaystyle a$};
\draw (390,150) node [anchor=north west][inner sep=0.75pt]   [align=left] {$\displaystyle c$};
\end{tikzpicture}
\caption{$X'(1+2)$}
\label{Skew-1}
\end{minipage}%
\begin{minipage}{.5\textwidth}
\centering
\tikzset{every picture/.style={line width=0.75pt}} 

\begin{tikzpicture}[x=0.75pt,y=0.75pt,yscale=-1,xscale=1]

\draw   (408,138.34) -- (383.59,197.27) -- (324.66,221.68) -- (265.73,197.27) -- (241.32,138.34) -- (265.73,79.41) -- (324.66,55) -- (383.59,79.41) -- cycle ;
\draw [color={rgb, 255:red, 255; green, 4; blue, 8 }  ,draw opacity=1 ][fill={rgb, 255:red, 240; green, 0; blue, 28 }  ,fill opacity=1 ] [dash pattern={on 0.84pt off 2.51pt}]  (324.66,58) -- (324.66,218.68) ;
\draw [shift={(324.66,221.68)}, rotate = 270] [fill={rgb, 255:red, 232; green, 1; blue, 29 }  ,fill opacity=1 ][line width=0.08]  [draw opacity=0] (8.93,-4.29) -- (0,0) -- (8.93,4.29) -- cycle    ;
\draw [shift={(324.66,55)}, rotate = 90] [fill={rgb, 255:red, 232; green, 1; blue, 29 }  ,fill opacity=1 ][line width=0.08]  [draw opacity=0] (8.93,-4.29) -- (0,0) -- (8.93,4.29) -- cycle    ;
\draw [color={rgb, 255:red, 232; green, 1; blue, 29 }  ,draw opacity=1 ][fill={rgb, 255:red, 232; green, 1; blue, 29 }  ,fill opacity=1 ] [dash pattern={on 0.84pt off 2.51pt}]  (244.32,138.34) -- (405,138.34) ;
\draw [shift={(408,138.34)}, rotate = 180] [fill={rgb, 255:red, 232; green, 1; blue, 29 }  ,fill opacity=1 ][line width=0.08]  [draw opacity=0] (8.93,-4.29) -- (0,0) -- (8.93,4.29) -- cycle    ;
\draw [shift={(241.32,138.34)}, rotate = 360] [fill={rgb, 255:red, 232; green, 1; blue, 29 }  ,fill opacity=1 ][line width=0.08]  [draw opacity=0] (8.93,-4.29) -- (0,0) -- (8.93,4.29) -- cycle    ;

\draw (320,225) node [anchor=north west][inner sep=0.75pt]   [align=left] {$\displaystyle a$};
\draw (320,40) node [anchor=north west][inner sep=0.75pt]   [align=left] {$\displaystyle a$};
\draw (385,195.27) node [anchor=north west][inner sep=0.75pt]   [align=left] {$\displaystyle b'$};
\draw (257,65) node [anchor=north west][inner sep=0.75pt]   [align=left] {$\displaystyle b$};
\draw (409,130) node [anchor=north west][inner sep=0.75pt]   [align=left] {$\displaystyle d$};
\draw (228,130) node [anchor=north west][inner sep=0.75pt]   [align=left] {$\displaystyle d$};
\draw (385,63) node [anchor=north west][inner sep=0.75pt]   [align=left] {$\displaystyle c$};
\draw (252,194) node [anchor=north west][inner sep=0.75pt]   [align=left] {$\displaystyle c'$};

\end{tikzpicture}

    \caption{$X'(2+4)$}
    \label{Skew-2}
\end{minipage}
\end{figure}

%% file: Examples.tex
\section{The Known Examples}
Before proving that $\TC(\al,1-\al)$, $Q_2(\frac{1}{3})$, and $Q_2(\frac{1}{3})^\times \oplus \GenG{\id}$ are the only examples (up to isomorphism) of $2$-generated axial algebras of Monster type with skew axet $X'(1+2)$, it would be useful to show that they satisfy the conditions we desire. 
\subsection{Three Dimensions}
Let $A=\TC(\al)$ with $\al\notin \{0,1, \frac{1}{2},-1\}$. In \cite{mcinroy2022axial} and \cite{franchi20212}, this example has already been stated. Let $x,y,z$ be the distinct axes of $\mcal{J}(\al)$-type and notice that $\{x,y,z\}$ form a basis for $A$. The element $\id=\frac{1}{\al+1}(x+y+z)$ is the identity element of $A$. Letting $w=\id-x$, one can see that $w$ is an axis of $\mcal{J}(1-\al)$-type and $\al\id=(-w+y+z)$. Looking at $\GenA{w,y}$, we have
\begin{eqnarray*}
wy &=& (\id-x)y = y-\frac{\al}{2}(x+y-z) = y-\frac{\al}{2}(-w+\id +y-z)\\
&=&y-\frac{\al}{2}(-w+y-z)-\frac{1}{2}(-w+y+z) = \frac{\al+1}{2}w +\frac{1-\al}{2}(y-z).
\end{eqnarray*}
Hence $z\in \GenA{w,y}$ and $x$ is too. Therefore $A=\GenA{w,y}$. As $A$ is now of $\mcal{M}(\al,1-\al)$-type with respect to the axes $w$ and $y$, let the Miyamoto involutions of $w$ and $y$ be $\tau_w$ and $\tau_y$ respectively. As $y$ is of $\mcal{J}(\al)$-type, $\tau_y$ is the identity. On the other hand, one can check that $x\in A_0(w)$ and $y-z\in A_{1-\al}(w)$. Further,
\[ y = \frac{\al}{2}x + \frac{\al+1}{2}w +\frac{1}{2}(y-z)\]
and so
\[ y^{\tau_w}=\frac{\al+1}{2}w+\frac{\al}{2}x -\frac{1}{2}(y-z) =z.\]
Therefore $A$ is skew with axet $X'(1+2)$. To avoid confusion about fusion laws, we denote the skew algebra by $\TC(\al,1-\al)$.

\begin{rem}
If one used the same method to $\TC(\frac{1}{2})$, then it would produce the algebra $S(1)$. By Remark 5.9 and Lemma 5.12 in \cite{mcinroy2021forbidden}, $S(1)$ has axet $X(6)$.
\end{rem}

\subsubsection{A Variation}
One may notice that for $\TC(-1)$, we have only $\jor{-1}$-axes and no identity. It would be absurd to try to make an axial algebra of Monster type out of it. We therefore use $\TC(2)$ to produce an axial algebra with fusion law $\mon{-1,2}$. Take $A=\TC(2)$ with $\{u,v,w\}$ being the three distinct axes of $\jor{2}$-type. There is an identity in $A$, $\id=\frac{1}{3}(u+v+w)$. Let $y:=\id-u$ and $z:=\id-v$ be $\jor{-1}$-axes. Let us look at $\GenA{w,y}$. We get
\[wy=(\id-u)w=w-uw=w-(u+w-v)=v-u\]
and
\[y(u-v)=(\id-u)(u-v)=-v+uv=-v+(u+v-w)=u-w.\]
Hence $v-w\in \GenA{w,y}$ as well as $v\in \GenA{w,y}$. Further, $u\in \GenA{w,y}$ and so $A=\GenA{w,y}$. These axes now satisfy the fusion law of $\mon{-1,2}$. We have that $\tu{y}$ is trivial as it is a $\jor{-1}$-axis while $\tu{w}$ is not. The reader can check:
\begin{itemize}
    \item $w \in A_1(w)$,
    \item $\id-w =-(y+z)\in A_0(w)$, and
    \item $y-z \in A_2(w)$.
\end{itemize}
We have
\[ y= \frac{1}{2}(y+z)+\frac{1}{2}(y-z),\]
and 
\[ y^{\tu{w}}=\frac{1}{2}(y+z)-\frac{1}{2}(y-z)=z.\]
Therefore $A$ is skew with axet $X'(1+2)$, and denote this algebra by $\TC(-1,2)$.
\begin{note}
We have $yz=(\id-u)(\id-v)=\id-u-v+(u+v)-w=\id-w=-y-z$ and so $\GenA{y,z}\cong \TC(-1)^\times$, the $2$-dimensional quotient of $\TC(-1)$. This will be important in the proof later.
\end{note}
\begin{rem}
Another way to construct $\TC(-1,2)$ would be by attaching a universal identity to $\TC(-1)^\times$ (similar to the construction of $Q_2(\frac{1}{3})^\times \oplus \GenG{\id}$).
\end{rem}

\subsection{Four Dimensions}
Let $\F$ have characteristic not equal to $3$ and let $Q_2(\frac{1}{3})$ be defined as in \cite[Section 5.3]{galt2021double} with $\eta=\frac{1}{3}$. The algebra has a basis of two single axes, $s_1$ and $s_2$, and two double axes, $d_1$ and $d_2$ with the multiplication defined in Table \ref{mult Q}.
\begin{table}[b]
\begin{center}
    \begin{tabular}{c|c|c|c|c}
         & $s_1$ & $s_2$ & $d_1$ &$d_2$   \\
         \hline
        $s_1$ & $s_1$ &$0$& $\frac{1}{3}s_1+\frac{1}{6}d_1-\frac{1}{6}d_2$& $\frac{1}{3}s_1-\frac{1}{6}d_1+\frac{1}{6}d_2$\\
        \hline
        $s_2$ & & $s_2$ & $\frac{1}{3}s_2+\frac{1}{6}d_1-\frac{1}{6}d_2$& $\frac{1}{3}s_2-\frac{1}{6}d_1+\frac{1}{6}d_2$\\
        \hline
        $d_1$ &  &  & $d_1$&$-\frac{1}{3}s_1-\frac{1}{3}s_2+\frac{1}{3}d_1+\frac{1}{3}d_2$ \\
        \hline
        $d_2$ &  &  &  & $d_2$
    \end{tabular}
       \caption{The multiplication table of $Q_2(\frac{1}{3})$}
       \label{mult Q}
       \end{center}
\end{table}
\subsubsection{In characteristic not equal to 5}
When the field has characteristic not equal to $5$, $A=Q_2(\frac{1}{3})$ is simple by Theorem 1.6 in \cite{galt2021double} and has identity element, $\id = \frac{3}{5}(s_1+s_2+d_1+d_2)$. Notice  $s_1$ is a $\jor{\frac{1}{3}}$-axis while $d_1$ is a $\mon{\frac{2}{3},\frac{1}{3}}$-axis. Let $t_1:=\id-d_1$ and $t_2:=\id -d_2$. Then $t_1$ and $t_2$ are $\mcal{M}(\frac{1}{3},\frac{2}{3})$-axes. Let us look at $\GenA{t_1,s_1}$. We have that 
\[ s_1t_1 = s_1(\id -d_1)=s_1-\frac{1}{3}s_1-\frac{1}{6}d_1+\frac{1}{6}d_2 = \frac{2}{3}s_1-\frac{1}{6}d_1+\frac{1}{6}d_2= \frac{2}{3}s_1+\frac{1}{6}t_1-\frac{1}{6}t_2.\]
Therefore $t_2\in \GenA{s_1,t_1}$. Notice that $-\frac{1}{3}\id = s_1+s_2-t_1-t_2$. We have
\begin{eqnarray*}t_1t_2 = (\id-d_1)(\id-d_2) &=& \id -d_1-d_2+d_1d_2\\
&=& \id -\frac{1}{3}s_1-\frac{1}{3}s_2-\frac{2}{3}d_1-\frac{2}{3}d_2\\
&=& \id -\frac{1}{3}s_1-\frac{1}{3}s_2+\frac{2}{3}(\id-d_1)+\frac{2}{3}(\id-d_2)-\frac{4}{3}\id\\
&=&-\frac{1}{3}\id -\frac{1}{3}s_1-\frac{1}{3}s_2+\frac{2}{3}t_1+\frac{2}{3}t_2\\
&=& s_1+s_2-t_1-t_2 -\frac{1}{3}s_1-\frac{1}{3}s_2+\frac{2}{3}t_1+\frac{2}{3}t_2\\
&=& \frac{2}{3}s_1+\frac{2}{3}s_2-\frac{1}{3}t_1-\frac{1}{3}t_2.
\end{eqnarray*}
Hence $s_2$ and $\id$ are in $\GenA{s_1,t_1}$. Therefore $d_1,d_2\in \GenA{s_1,t_1}$ and so $A=\GenA{s_1,t_1}$ and has fusion law $\mon{\frac{1}{3},\frac{2}{3}}$. As $s_1$ is a $\jor{\frac{1}{3}}$-axis, $\tu{{s_1}}$ is the identity map.  The reader can check that:
\begin{itemize}
    \item $t_1 = \frac{1}{5}(3s_1+3s_2-2d_1+3d_2)\in A_1(t_1)$,
    \item $d_1 \in A_0(t_1)$, 
    \item $s_1+s_2-d_2 \in A_\frac{1}{3}(t_1)$, and
    \item $s_1-s_2 \in  A_\frac{2}{3}(t_1)$.
\end{itemize}
Hence
\[s_1 = \frac{5}{12}t_1+\frac{1}{6}d_1+\frac{1}{4}(s_1+s_2-d_2)+\frac{1}{2}(s_1-s_2).\]
Therefore
\[ s_1^{\tu{t_1}} = \frac{5}{12}t_1+\frac{1}{6}d_1+\frac{1}{4}(s_1+s_2-d_2)-\frac{1}{2}(s_1-s_2)=s_2.\]
Hence $A$ is skew with axet $X'(1+2)$. To avoid confusion with fusion laws, we denote this skew algebra by $Q_2(\frac{1}{3},\frac{2}{3})$. 
\begin{note}
One may think we can make any $Q_2(\eta)$, $\eta\notin \{-\frac{1}{2},\frac{1}{3}\}$, a skew axial algebra of Monster type with this method. This is partly possible if $\GenA{t_1,s_1}=Q_2(\eta)$. However the algebra would be skew but not of Monster type. The fusion law would have five elements and be an extension of the Monster fusion law.
\end{note}
\subsubsection{In characteristic 5}
Suppose now that the field has characteristic $5$. Then $\frac{1}{3}=-\frac{1}{2}$ and $Q_2(\frac{1}{3})$ is not simple by Theorem 1.6 in \cite{galt2021double}. This algebra has an annihilating element rather than an identity and so it is impossible with the technique used so far. We have that $I=\GenG{s_1+s_2+d_1+d_2}$ is the radical of $Q_2(\frac{1}{3})$. The quotient is denoted by $Q_2(\frac{1}{3})^\times$ and it is spanned by axes $\{x,y,z\}$ with their multiplication stated in Table \ref{mult Qx} (without the final row and column).
\begin{table}[b]
\begin{center}
    \begin{tabular}{c|c|c|c|c}
         & $x$ & $y$ & $z$ & $\id$   \\
         \hline
        $x$ & $x$ &$0$& $3x+y+2z$ & $x$\\
        \hline
        $y$ & & $y$ & $x+3y+2z$ & $y$\\
        \hline
        $z$ &  &  & $z$ & $z$\\
        \hline 
        $\id$ & & & & $\id$
    \end{tabular}
       \caption{The multiplication table of $Q_2(\frac{1}{3})^\times \oplus \GenG{\id}$}
       \label{mult Qx}
       \end{center}
\end{table}
The reader can check that $Q_2(\frac{1}{3})^\times =\GenA{x,z}$ with $x$ a $\jor{\frac{1}{3}}$-axis and $z$ a $\mon{\frac{2}{3},\frac{1}{3}}$-axis and so cannot be skew. In fact, this algebra has an axet of type $X(4)$. 

Let $A=Q_2(\frac{1}{3})^\times \oplus \GenG{\id}$ to be a $4$-dimensional algebra, where we add a universal identity element, $\id$, and multiplication is stated in Table \ref{mult Qx}. Let $w:=\id -z$ and let us look at $\GenA{x,w}$. Notice that
\[ wx=x(\id-z)=x - (3x+y+2z)=-2x-y-2z.\]
Therefore $y+2z\in \GenA{x,w}$ and we have
\[ w(y+2z)=wy=(\id-z)y=y-(x+3y+2z)=-x-2y-2z.\]
Hence $y+z\in \GenA{x,w}$. Moreover, $y,z \in \GenA{x,w}$ and so is $\id$. Therefore $A=\GenA{w,x}$. Showing that $x$ and $w$ are axes is a straightforward task as $x$ and $z$ are axes of the subalgebra $Q_2(\frac{1}{3})^\times$. We have that $x$ is a $\jor{\frac{1}{3}}$-axis and $w$ is a $\mon{\frac{1}{3},\frac{2}{3}}$-axis. Therefore $\tu{x}$ is trivial and $\tu{w}$ is not. For $w$, we have the following eigenvectors:
\begin{itemize}
    \item $w\in A_1(w)$,
    \item $z \in A_0(w)$,
    \item $x+y+3z \in A_{\frac{1}{3}}(w)$, and
    \item $x-y \in A_{\frac{2}{3}}(w)$. 
\end{itemize}
We have
\[ x= z+\frac{1}{2}(x+y+3z)+\frac{1}{2}(x-y)\]
and so
\[ x^{\tu{w}}=z+\frac{1}{2}(x+y+3z)-\frac{1}{2}(x-y)=y.\]
Therefore $A$ is skew with axet $X'(1+2)$. Further, $A\not\cong Q_2(\frac{1}{3})$ as $A$ has an identity element while $Q_2(\frac{1}{3})$ has not, in characteristic $5$. 

\begin{note}
One may think we can apply this method to $Q_2(-\frac{1}{2})^\times$ in characteristic not equal to $5$ to produce new algebras. This can be done however we run in the same problem as $Q_2(\eta)$ where the fusion law is an extension of the Monster fusion law and one needs to check that it is $2$-generated.
\end{note}

%% file: Construction.tex
\section{The Foundations}
This section will be working closely with \cite[Section 4]{franchi20211} and \cite[Section 3]{rehren2017generalised}, which will be cited accordingly. We will state some of the results in those papers however we will omit the proofs.

Fix $(A,X)$ to be an $\mon{\al,\bt}$-axial algebra with $X=\{a_0,a_1\}$. For $i\in\{0,1\}$, let $\tu{i}$ be the Miyamoto involution associated to $a_i$. Let $\GenG{X}$ be isomorphic to $X'(1+2)$ and without loss of generality, let $\tu{0}$  be a non-trivial involution and $\tu{1}$ be equal to the identity. Set $a_{2i}:=a_0^{\tu{0}^i}$ and $a_{2i+1}:=a_1^{\tu{0}^i}$ for $i\in \Z$. Hence, 
\begin{equation}\label{axes rel}
    a_{2i}=a_0,\;\; a_{4i+1}=a_1,\;\; a_{4i-1}=a_{-1}\; \text{for} \;i \in \Z.
\end{equation} 
For $i\in \Z$ and $m\in \N$, define 
\begin{equation}\label{s rel}
 s_{i,m}:=a_ia_{i+m}-\bt(a_i+a_{i+m}).   
\end{equation} Let $f$ be the semi-automorphism of $A$ such that $a_i^f=a_{1-i}$ for all $i\in \Z$. This map is well-defined by Corollary 3.8 in \cite{franchi20211}. This $f$ is called the \emph{flip} and we do not assume it is an automorphism of $A$. Notice that $f$ has order $2$, $s_{0,1}^f=s_{0,1}$ and $s_{0,2}^f=s_{1,2}$. Using Definition \ref{Proj}, we define \[\lm_1:=\lm_{a_0}(a_1), \; \; \lmf_1:=\lm_{a_1}(a_0), \; \; \lm_2:=\lm_{a_0}(a_2)=1, \; \;\text{and }
\lmf_2:=\lm_{a_1}(a_{-1}).\]
To ease notation, we define the following constants:
\[ \gm:=\bt-\lm_1,\; \ep:=(1-\al)\lm_1-\bt, \; \dt:=(1-\al)\lm_1+\bt(\al-\bt-1).\]
In a similar fashion,
\[ \gm^f:=\bt-\lmf_1, \; \ep^f:=(1-\al)\lmf_1-\bt, \; \dt^f:=(1-\al)\lmf_1+\bt(\al-\bt-1).\]
They have the following relations:
\[
(\al-1)\gm=\ep+\al\bt=\dt+\bt^2\;\text{ and } \;(\al-1)\gm^f=\ep^f+\al\bt=\dt^f+\bt^2.
\]

\begin{lem}\label{axes sigma}
We have 
\begin{eqnarray*}
a_0s_{0,1}&=&(\al-\bt)s_{0,1}+\dt a_0 +\frac{1}{2}\bt(\al-\bt)(a_1+a_{-1})\\
a_1 s_{0,1}&=& (\al-\bt)s_{0,1}+ \dt^f a_1+\bt(\al-\bt) a_0 \\
a_{-1} s_{0,1} &=& (\al-\bt)s_{0,1} + \dt^f a_{-1}+\bt(\al-\bt) a_0.
\end{eqnarray*}
\proof The first two equalities are part of Lemma 3.1 in \cite{rehren2017generalised}. 
By the second equality, 
\[(a_1s_{0,1})^{\tu{0}} = [(\al-\bt)s_{0,1}+ \dt^f a_{1}+\bt(\al-\bt) a_0 ]^{\tu{0}}\]
if and only if
\[ a_{-1}s_{0,1} = (\al-\bt)s_{0,1}+ \dt^f a_{-1}+\bt(\al-\bt) a_0 \]
due to $\tu{0}$ being an automorphism of $A$. \qed
\end{lem}
\begin{lem}\label{sigma2}
We have that $s_{1,2}$ is in the span of $\{a_0,a_1,a_{-1}, s_{0,1}\}$.
\proof Applying $f$ to the first equality of Lemma 4.7 in \cite{franchi20211}, we get:
\begin{eqnarray*}
(\al-2\bt)a_1s_{0,2}&=&\bt^2(\al-\bt)(a_3+a_{-1})+2\bt(\al-\bt)s_{1,2}\\
&+&\left[-2\al\bt\lmf_1+2\bt(1-\al)\lm_1\right.\\
&+&\left.\frac{\bt}{2}(4\al^2-2\al\bt+4\bt^2-\al-2\bt)\right](a_0+a_2)\\
&+&\frac{1}{(\al-\bt)}\left[(6\al^2-8\al\bt-2\al+4\bt)(\lmf_1)^2\right.\\
&+&\left. 2\al(\al-1)\lm_1\lmf_1+2\al(-2\al-2\bt+1)(\al-\bt)\lmf_1\right.\\
&-&\left.4\bt(\al-1)(\al-\bt)\lm_1-\al\bt(\al-\bt)\lmf_2\right.\\
&+&\left.2\bt(2\al^2+\bt^2-\al)(\al-\bt)\right]a_1\\
&+&\left[-4\al\lmf_1-4(\al-1)\lm_1\right.\\
&+&\left.(4\al^2-2\al\bt+4\bt^2-\al-2\bt)\right]s_{0,1}.
\end{eqnarray*}
Applying relations in $(\ref{axes rel})$ and the fact that $s_{0,2}=a_0a_2-\bt(a_0+a_2)=(1-2\bt)a_0$ we get
\begin{eqnarray*}
(\al-2\bt)(1-2\bt)a_1a_0&=&2\bt^2(\al-\bt)a_{-1}+2\bt(\al-\bt)s_{1,2}\\
&+&2\left[-2\al\bt\lmf_1+2\bt(1-\al)\lm_1\right.\\
&+&\left.\frac{\bt}{2}(4\al^2-2\al\bt+4\bt^2-\al-2\bt)\right]a_0\\
&+&\frac{1}{(\al-\bt)}\left[(6\al^2-8\al\bt-2\al+4\bt)(\lmf_1)^2\right.\\
&+&\left. 2\al(\al-1)\lm_1\lmf_1+2\al(-2\al-2\bt+1)(\al-\bt)\lmf_1\right.\\
&-&\left.4\bt(\al-1)(\al-\bt)\lm_1-\al\bt(\al-\bt)\lmf_2\right.\\
&+&\left.2\bt(2\al^2+\bt^2-\al)(\al-\bt)\right]a_1\\
&+&\left[-4\al\lmf_1-4(\al-1)\lm_1\right.\\
&+&\left.(4\al^2-2\al\bt+4\bt^2-\al-2\bt)\right]s_{0,1}.\\
\end{eqnarray*}
Equivalently,
\begin{eqnarray*}
-2\bt(\al-\bt)s_{1,2}&=&2\bt^2(\al-\bt)a_{-1}\\
&+&2\left[-2\al\bt\lmf_1+2\bt(1-\al)\lm_1\right.\\
&+&\left.\frac{\bt}{2}(4\al^2-2\al\bt+4\bt^2-\al-2\bt)\right]a_0\\
&+&\frac{1}{(\al-\bt)}\left[(6\al^2-8\al\bt-2\al+4\bt)(\lmf_1)^2\right.\\
&+&\left.2\al(\al-1)\lm_1\lmf_1+2\al(-2\al-2\bt+1)(\al-\bt)\lmf_1\right.\\
&-&\left.4\bt(\al-1)(\al-\bt)\lm_1-\al\bt(\al-\bt)\lmf_2\right.\\
&+&\left.2\bt(2\al^2+\bt^2-\al)(\al-\bt)\right]a_1\\
&+&\left[-4\al\lmf_1-4(\al-1)\lm_1\right.\\
&+&\left.(4\al^2-2\al\bt+4\bt^2-\al-2\bt)\right]s_{0,1}\\
&-&(\al-2\bt)(1-2\bt)\left[s_{0,1}+\bt(a_0+a_1)\right].
\end{eqnarray*}
Therefore
\[ s_{1,2}=Pa_0+Qa_1+Ra_{-1}+Ss_{0,1}\]
where 
\begin{eqnarray*}
P&:=& \frac{1}{(\al-\bt)}\left[2(\al-1)\lm_1+2\al\lmf_1+\al(1-2\al)\right],\\
Q&:=&-\frac{1}{2\bt(\al-\bt)^2}\left[(6\al^2-8\al\bt-2\al+4\bt)(\lmf_1)^2+2\al(\al-1)\lm_1\lmf_1\right.\\
&+&\left.2\al(-2\al-2\bt+1)(\al-\bt)\lmf_1-4\bt(\al-1)(\al-\bt)\lm_1\right.\\
&-&\left.\al\bt(\al-\bt)\lmf_2+2\bt(2\al^2+\bt^2-\al)(\al-\bt)\right.\\
&-&\left.\bt(\al-\bt)(\al-2\bt)(1-2\bt)\right],\\
R&:=&-\bt, \text{ and}\\
S&:=&\frac{P}{\bt}.
\end{eqnarray*}
\qed
\end{lem}
\begin{lem}\label{mult lemma}
With the above notation, $Q=R$ and $a_1a_{-1}=P(a_0+\frac{1}{\bt}\sg)$.
\end{lem}
\proof
Since $s_{1,2}$ is invariant under $\tu{0}$, notice that 
\[0=(s_{1,2}-s_{1,2}^{\tu{0}})=(Q-R)(a_1-a_{-1})\]
To avoid a contradiction, $Q=R$ and so
\[ s_{1,2}=P\left(a_0+\frac{1}{\bt}s_{0,1}\right)-\bt(a_1+a_{-1}).\]
Hence
\[a_1a_{-1}= s_{1,2}+\bt(a_1+a_{-1})=P\left(a_0+\frac{1}{\bt}s_{0,1}\right).\] \qed
\begin{prop}\label{span prop}
We have $A$ is linearly spanned by $\mcal{B}=\{a_{-1},a_0,a_1,s_{0,1}\}$ and so is at most $4$-dimensional.
\proof From above, multiplication of $\mcal{B}$ has been described and shown to be in $\Span{\mcal{B}}$ except for $s_{0,1}^2$. If $\al\neq 2\bt$,  $s_{0,1}^2 \in \Span{\mcal{B}}$ by Lemma 4.7 in \cite{franchi20211}. 
If $\al = 2\bt$, by Lemma 3.5 in \cite{franchi20212}, $s_{0,1}^2$ is computed. As $s_{0,3}=s_{0,1}$ by Equation $(\ref{s rel})$, $s_{0,1}^2\in \Span{\mcal{B}}$.  \qed
\end{prop}

We will now state  Theorem 4.1.1 in \cite{rehren2015axial} where $R$ is a field. This is a extremely useful result that will be used multiple times in our reasoning. 
\begin{prop}[Rehren]\label{Rehren Thm}
Let $\al,\bt\notin \{0,1\}$ be distinct values. Suppose that $V=\GenA{p,q}$ be an axial algebra such that $p$ is a $\mcal{J}(\al)$-axis and $q$ is a $\mcal{J}(\bt)$-axis. Then $V\cong \TB$ or $V\cong \TC(\al,1-\al)$.
\end{prop}

\begin{rem}
In Rehren's proof, he only shows algebra isomorphisms of $V\cong \TC(\al)$ or $V\cong \TC(2)$ if $\al=-1$. However, for $\al\neq-1$, it is clear that $V\cong \TC(\al,1-\al)$ due $(V,\{p,q\})$ satisfying a fusion law $\mon{\al,1-\al}$. Further, for $\al=-1$, it is easy to show that $V\cong \TC(-1,2)$. Again, this is clear since $(V,\{p,q\})$ satisfies a fusion law $\mon{-1,2}$.  
\end{rem}

%% file: SkewRelations.tex
\section{Skew Relations}
To ease notation, set $a:=a_0$, $b:=a_1$, $c:=a_{-1}$, and $\sg:=s_{0,1}$. From Lemma \ref{axes sigma} and \ref{sigma2}, we have
\begin{eqnarray*}
a\sg&=&(\al-\bt)\sg+\dt a+\frac{1}{2}\bt(\al-\bt)(b+c),\\
b\sg&=& (\al-\bt)\sg+\bt(\al-\bt)a + \dt^f b,\\
c\sg&=& (\al-\bt)\sg+\bt(\al-\bt)a + \dt^f c,\\
 bc &=& P\left(a+\frac{1}{\bt}\sg\right).
\end{eqnarray*}
The last thing  to complete a multiplication table of $A$ is $\sg^2$. These expressions are in \cite{franchi20211} and \cite{franchi20212} for $\al\neq 2\bt$ and $\al=2\bt$ respectively. In this paper, we do not need the exact expression for $\sg^2$ and we simply write
\[ \sg^2 = \zt a + \te(b+c)+\kp \sg\]
with $\zt, \te, \kp \in \F$.
Notice that the $b$ and $c$ term are equal due to $\sg^2$ being invariant under $\tu{0}$ and Table \ref{mult table} gives us the complete multiplication of $A$.
\begin{table}[b]
\centering
\scalebox{0.9}{
    \begin{tabular}{c|c|c|c|c}
         & $a$ & $b$ & $c$ &$\sg$   \\
         \hline
        $a$ & $a$& $\bt a + \bt b + \sg$& $\bt a + \bt c+ \sg$& $\dt a+\frac{1}{2}\bt(\al-\bt)(b+c)+(\al-\bt)\sg$\\
        \hline
        $b$ & & $b$ & $P(a+\frac{1}{\bt}\sg)$ &$\bt(\al-\bt)a + \dt^f b+(\al-\bt)\sg$\\
        \hline
        $c$ &  &  & $c$& $\bt(\al-\bt)a + \dt^f c+(\al-\bt)\sg$\\
        \hline
        $\sg$ &  &  &  & $\zt a + \te(b+c)+\kp \sg$
    \end{tabular}}
\caption{The multiplication table of $A$}
\label{mult table}
\end{table}

The rest of this section will be producing results that will be used in the proof of Theorem \ref{Theorem}. First, we can calculate eigenspaces of $\text{ad}_a$ and $\text{ad}_b$.

\begin{lem}\label{eigen a}
The eigenspaces of $\text{ad}_a$ are the following:
\begin{enumerate}
    \item[$1$.] $\GenG{a}=A_1(a)$, 
    \item[$2$.] $\GenG{\ep a+\frac{1}{2}(\al-\bt)(b+c)-\sg}= A_0(a)$,
    \item[$3$.] $\GenG{\gm a +\frac{1}{2}\bt(b+c)+\sg}=A_\al(a)$, and
    \item[$4$.] $\GenG{b-c}= A_\bt(a)$.
\end{enumerate}
\proof See Lemma 4.4 in \cite{franchi20211}. \qed
\end{lem}
\begin{lem}\label{eigen b}
The eigenspaces of $\text{ad}_b$ are the following:
\begin{itemize}
    \item[$1$.] $\GenG{b}= A_1(b)$,
    \item[$2$.] $\GenG{-\frac{P}{\bt}a+Pb+c, (\al-\bt)a +\ep^fb-\sg}=A_0(b)$, and
    \item[$3$.] $\GenG{\bt a +\gm^fb+\sg}=A_\al(b)$.
\end{itemize}
\proof Using Table \ref{mult table}, the reader can check the generators are eigenvectors. As $A$ is at most $4$-dimensional, there are no more possible eigenvectors (up to linear combination). \qed 
\end{lem}
We can find certain relations from the conditions imposed by Lemma \ref{Seress}. To avoid focusing on how these relations are calculated, the reasoning is in the appendix for the reader to look at if they wish. 
\begin{lem}\label{skew rel}
We have that the following equations must hold:
\begin{equation}\label{proof1}
   \begin{split}
\lmf_2=-\frac{P}{\bt}\gm^f
\end{split} 
\end{equation} 
\begin{equation}\label{proof2}
    \begin{split}
\bt\dt =\frac{1}{2}\bt(\al-\bt)-\bt^2(\al-\bt)-(\al-2\bt)\dt^f.  
\end{split}
\end{equation} 
and
\begin{equation}\label{proof3}
\begin{split}
\frac{1}{2}(1-\bt)P=(\al-1)\gm^f.
\end{split}
\end{equation}

\end{lem}
\proof See the Appendix. \qed

%% file: TheoremProof.tex
\section{Proof of Theorem \ref{Theorem}}

Let $(A,X)$ be an $\mon{\al,\bt}$-axial algebra satisfying the hypothesis of Theorem \ref{Theorem} with $X=\{a,b\}$. By Proposition \ref{span prop}, we may assume that $A$ is spanned by $\{a,b,c,\sg\}$ and multiplication is defined in Table \ref{mult table}. We assume $b$ is a $\jor{\al}$-axis with $c=b^{\tu{a}}$. We set 
$U:=\GenA{b,c}$ which is a  $\jor{\al}$-axial algebra. The proof is split into two parts: $P=0$ or $P\neq 0$. 

\subsection{Orthogonal axes}\label{P zero}
For this section, we will assume that $P=0$.
By Table \ref{mult table}, $bc=0$; that is $U\cong \TB$. Hence $d:=b+c$ is a $\jor{2\al}$-axis in $A$. Notice that $d$ is not primitive as $A_1(d)=\GenG{b,c}$. Let $F:=A_{\{0,1,\al\}}(a)$ be the fixed subalgebra of $\tu{a}$. We have that $a,d,\sg \in F$ and $b-c\notin F$ which implies $\text{dim}F\leq\text{dim}A - 1\leq 3$.

\begin{lem}
We have $\al=\frac{1}{3}$, $F\cong \TC(\frac{1}{3},\frac{2}{3})$ or $F\cong \TC(\frac{2}{3},\frac{1}{3})$, and $A$ has dimension $4$.
\proof  One can show $F=\GenA{a,d}$ is a primitive axial algebra. We know $a$ satisfies $\jor{\al}$ fusion law and $d$ satisfies the $\jor{2\al}$ fusion law. As $\al\neq 2\al$, then $F\cong \TB$, $F\cong \TC(\frac{1}{3},\frac{2}{3})$, or $F\cong \TC(\frac{2}{3},\frac{1}{3})$ by Proposition \ref{Rehren Thm}. Let us assume $F\cong \TB$. Then $ad=0$ thus $\sg =-\bt a -\frac{\bt}{2}(b+c)$. Therefore the $\al$-eigenvector of $\text{ad}_a$ by Lemma \ref{eigen a} is now
\[ \gm a +\frac{1}{2}\bt(b+c)+\sg = (\gm-\bt)a=-\lm_1 a.\]
To avoid contradiction, we must have $\lm_1=0$ and $a$ is a $\jor{\bt}$-axis in $A$. By Proposition \ref{Rehren Thm}, $A\cong \TC(\al,1-\al)$ and $\id$ exists. Since $bc=0$, we must have $c=\id -b$ but that would make $c$ a $\jor{\bt}$-axis in $A$, which would contradict our axet. Therefore $\al=\frac{1}{3}$. We have $F\cong \TC(\frac{1}{3},\frac{2}{3})$ or $F\cong \TC(\frac{2}{3},\frac{1}{3})$. Since $F$ is $3$-dimensional, $A$ has to be $4$-dimensional. \qed
\end{lem}

\begin{lem}
We have that $A$ has fusion law $\mon{\frac{1}{3}, \frac{2}{3}}$ and multiplication is shown in Table $\ref{mult 4dim}$.

\proof As $P=0$, we have that $\gm^f=0$, moreover $\lmf_1=\bt$, by Equation (\ref{proof3}). Substituting into Equation (\ref{proof2}), we have
\[\bt \dt = \frac{1}{2}\bt(\al-\bt)-\bt^2(\al-\bt)+\bt^2(\al-2\bt)=\bt\left(\frac{1}{2}(\al-\bt)-\bt^2\right).\]
With $\al=\frac{1}{3}$, we get $\lm_1=\frac{1}{4}(\bt+1)$. Further, with $P=0$, we get
\[ 0=-\frac{4}{3}\lm_1+\frac{2}{3}\bt+\frac{1}{9}=-\frac{1}{3}(\bt+1)+\frac{2}{3}\bt+\frac{1}{9}=\frac{1}{3}\bt-\frac{2}{9}.\]
Therefore $\bt=\frac{2}{3}$.

With $F\cong \TC(\frac{1}{3},\frac{2}{3})$ or $F\cong \TC(\frac{2}{3},\frac{1}{3})$, there are exactly three $\mcal{J}(\frac{1}{3})$-axes in $F$, which are $a$ and two other elements, $e$ and $f$. Without loss of generality, let $e=\id -d$ where $\id=\frac{3}{4}(a+e+f)=3(a-d+f)$. Therefore $\{a,d,f,b\}$ is a basis for $A$ and by the multiplication of $F$, we have
\begin{eqnarray*}
ad&=&a(\id-e)=a-ae=a-\frac{1}{6}(a+e-f)=\frac{5}{6}a+\frac{1}{6}f-\frac{1}{6}(\id-d)\\
&=&\frac{5}{6}a+\frac{1}{6}f+\frac{1}{6}d-\frac{1}{2}(a-d+f)= \frac{1}{3}a +\frac{2}{3}d-\frac{1}{3}f,
\end{eqnarray*}
\begin{eqnarray*}
fd&=&f(\id-e)=f-fe=f-\frac{1}{6}(f+e-a)=\frac{5}{6}f+\frac{1}{6}a-\frac{1}{6}(\id-d)\\
&=&\frac{5}{6}f+\frac{1}{6}a+\frac{1}{6}d-\frac{1}{2}(a-d+f)= \frac{1}{3}f +\frac{2}{3}d-\frac{1}{3}a,
\end{eqnarray*}
and
\begin{eqnarray*}
af&=&\frac{1}{6}(a+f-e)=\frac{1}{6}(a+f+d-\id)\\
&=&\frac{1}{6}(a+f+d-3(a-d+f))= -\frac{1}{3}a+\frac{2}{3}d-\frac{1}{3}f. 
\end{eqnarray*}
Further,
\[ bd = b(b+c)=b+bc=b.\]
We have $ad =a(b+c)= 2\sg +\frac{4}{3} a +\frac{2}{3} d$, by the multiplication table of $A$, and so
\[ \sg =-\frac{1}{2}a -\frac{1}{6}f.\]
Hence 
\[ ab = \sg +\frac{2}{3} a +\frac{2}{3} b = \frac{1}{6}a-\frac{1}{6}f+\frac{2}{3} b.\] 
As $\lmf_1=\bt$, notice that $\frac{1}{6}a-\frac{2}{9} b +\frac{1}{6}f \in A_0(b)$ by Lemma \ref{eigen b}. Thus
\[
0 = b\left(a-\frac{4}{3}b +f\right) = \frac{1}{6}a -\frac{1}{6}f+\frac{2}{3} b -\frac{4}{3} b +bf\\
\]
which is equivalent to
\[
bf =-\frac{1}{6}a +\frac{1}{6}f+\frac{2}{3}b.
\]
Using $c=d-b$, we get the multiplication in Table \ref{mult 4dim}.\qed
\end{lem}
\begin{table}[ht]
\begin{center}
    \begin{tabular}{c|c|c|c|c}
         & $b$ & $c$ & $a$ &$f$   \\
         \hline
        $b$ & $b$ &$0$& $\frac{2}{3}b+\frac{1}{6}a-\frac{1}{6}f$& $\frac{2}{3}b-\frac{1}{6}a+\frac{1}{6}f$\\
        \hline
        $c$ & & $c$ & $\frac{2}{3}c+\frac{1}{6}a-\frac{1}{6}f$& $\frac{2}{3}c-\frac{1}{6}a+\frac{1}{6}f$\\
        \hline
        $a$ &  &  & $a$&$\frac{2}{3}b+\frac{2}{3}c-\frac{1}{3}a-\frac{1}{3}f$ \\
        \hline
        $f$ &  &  &  & $f$
    \end{tabular}
       \caption{The multiplication table of $A$ when $P=0$}
       \label{mult 4dim}
       \end{center}
\end{table}
\begin{lem}
Let $\F$ have characteristic not equal to $5$. Then $A\cong Q_2(\frac{1}{3},\frac{2}{3})$.
\proof We define the map $\varphi: A \rightarrow Q_2(\frac{1}{3},\frac{2}{3})$ with $\varphi(a)=t_1$, $\varphi(b)=s_1$, $\varphi(c)=s_2$, and $\varphi(f)=t_2$. The reader can check that this is an isomorphism. \qed
\end{lem}
\begin{lem}
Let $\F$ have characteristic equal to $5$. Then $A\cong Q_2(\frac{1}{3})^\times \oplus \GenG{\id}$. 
\proof We have the map $\phi: A \rightarrow Q_2(\frac{1}{3})^\times\oplus \GenG{\id}$ with $\phi(a)=\id-z$, $\phi(b)=x$, $\phi(c)=y$, and $\phi(\id)=\id$. This is a straightforward check that $\phi$ is an isomorphism. \qed
\end{lem}

\subsection{Non-Orthogonal Case}\label{P not zero}
We will now assume $P\neq 0$. As $U$ is a axial algebra of $\jor{\al}$-type, $U$ is either $\TB$, $S(2)^\circ$, $\TC(-1)^\times$ or $3$-dimensional. We should note that $S(2)^\circ$ is $\text{Cl}^0(\F^2,b)$ in \cite{hall2015primitive}. We will look at each possibility separately to complete our proof. We recommend the reader to \cite[Section 3]{hall2015primitive} and  \cite[Section 5]{mcinroy2021forbidden} for more information on $2$-generated axial algebras of Jordan type.

\begin{lem}
We have $U\not\cong \TB$. 
\proof Let $U\cong \TB$. As $P\neq 0$ then $\sg = -\bt a$. Then $A$ is at most $3$-dimensional. Notice that $ab = \bt b$. Therefore $b\in A_\bt(a)$ and we have that $c=b^{\tu{a}}=-b$. This is a contradiction as $c^2=(-b)^2=b\neq c$. \qed
\end{lem}
\begin{lem}\label{S2}
We have $U\not\cong S(2)^\circ$.
\end{lem}
\proof Let $U\cong S(2)^\circ$ and $bc=\frac{1}{2}(b+c)$. Further $bc=P(a+\frac{1}{\bt}\sg)$. Since $P\neq 0$, we can express $\sg$ in terms of $a$, $b$ and $c$. Thus
\[ \sg = \frac{\bt}{2P}(b+c)-\bt a.\]
Note that
\[ab =\frac{\bt}{2P}(b+c)+\bt b\;\text{ and }\;ac = \frac{\bt}{2P}(b+c)+\bt c.\]
Hence $(b+c)$ is an eigenvector of $\text{ad}_a$ with eigenvalue $\mu=\frac{\bt}{P}+\bt$. By $a$ being primitive and $\bt\neq 0$, $\mu\notin \{ 1, \bt\}$. Hence $\mu \in \{0,\frac{1}{2}\}$.

Suppose $\mu=\frac{1}{2}$ and so $(b+c)\in A_\frac{1}{2}(a)$. By the Monster fusion law, $(b+c)^2 \in A_{\{0,1\}}(a)$. However
\[ (b+c)^2=b+c+2bc=2(b+c)\in A_\frac{1}{2}(a).\]
Thus $b+c=0$ which contradicts our choice of $U$. 

Therefore $\mu=0$ and $a$ is a $\jor{\bt}$-axis in $A$. As $b$ is a $\jor{\frac{1}{2}}$-axis and by Proposition \ref{Rehren Thm}, $\frac{1}{2}+\bt=1$. Thus $\bt=\frac{1}{2}=\al$ which contradicts the fusion law. \qed
\begin{lem}
If $U\cong \TC(-1)^\times$, then $A\cong \TC(-1,2)$.
\proof Suppose $U\cong \TC(-1)^\times$. We will use a similar method to Lemma \ref{S2}. We have $bc= -(b+c)$ and $bc=P(a+\frac{1}{\bt}\sg)$ thus
\[ \sg = -\bt a -\frac{\bt}{P}(b+c)\] 
and $A$ is at most 3-dimensional. Note that
\[ ab = -\frac{\bt}{P}(b+c)+\bt b\; \text{ and }\; ac = -\frac{\bt}{P}(b+c)+\bt c.
\]
Hence $(b+c)$ is an eigenvector of $\text{ad}_a$ with eigenvector $\mu=-\frac{2\bt}{P}+\bt$. Note that $\mu\notin \{1,\bt\}$ as it would contradict $a$ being primitive and $\bt\neq 0$ respectively. Thus $\mu \in \{0,-1\}$.

Suppose $\mu=-1$ moreover $(b+c)\in A_{-1}(a)$. By the Monster fusion law, $(b+c)^2\in A_{\{1,0\}}(a)$. However
\[ (b+c)^2=b+c+2bc=-(b+c)\in A_{-1}(a).\]
Thus $b+c=0$ which contradicts $U$.

Therefore $\mu=0$ and $a$ is $\jor{\bt}$-axis and by Proposition \ref{Rehren Thm}, we get $\bt=2$ and $A\cong \TC(-1,2)$. \qed
\end{lem}
If $U$ is $3$-dimensional, then $U$ is either $\TC(\al)$, $S(\dt)$ with $\dt\neq 2$ (denoted by $\text{Cl}^J(\F^2,b)$ in \cite{hall2015primitive}) or $\widehat{S}(2)^\circ$ (denoted by $\text{Cl}^{00}(\F^2,b)$ in \cite{hall2015primitive}).  

\begin{lem}
Let $U$ be $3$-dimensional, then $A=U$. Further, $U\cong \TC(\al)$ and $A\cong \TC(\al,1-\al)$ for $\al\neq -1$.
\proof Assume for a contradiction that $A\neq U$. Therefore $A$ is 4-dimensional with basis $\{a,b,c,\sg\}$. By Lemma \ref{eigen b}, we have
\begin{eqnarray*}
-\frac{P}{\bt}a +P b +c \in A_0(b)\text{ and } \bt a +\gm^f b +\sg \in A_\al(b).
\end{eqnarray*}
By the fusion law, the product of these two vectors will be in $A_\al(b)=\GenG{\bt a +\gm^fb +\sg}$ moreover the $c$ component is zero. We have
\begin{eqnarray*}
\left[-\frac{P}{\bt}a +P b +c \right]\left[\bt a +\gm^f b +\sg\right] &=& -Pa -\frac{P}{\bt}\gm^fab -\frac{P}{\bt}a\sg +P\bt ab\\
&+&P\gm^f b +Pb\sg +\bt ac +\gm^fbc +c\sg\\
&=&[...]a + [...]b+ [...]\sg\\
&+& \left[-\frac{1}{2}(\al-\bt)P+\bt^2+\dt^f\right]c.
\end{eqnarray*}
Thus
\[ \frac{1}{2}(\al-\bt)P =(\al-1)\gm^f.\]
However by Equation (\ref{proof3}), we get
\[\frac{1}{2}(\al-\bt)P =(\al-1)\gm^f  = \frac{1}{2}(1-\bt)P\]
which implies that $\al=1$. A contradiction and so $U=A$. 

The possible options for $U$ (and $A$) are $S(\dt)$ with $\dt\neq 2$, $\widehat{S}(2)^\circ$ and $\TC(\al)$. Let us look at each algebra separately.

Suppose $U\cong S(\dt)$. This is the same algebra in Section  3.5 of \cite{hall2015primitive}. It is shown that all non-trivial idempotents in $U$ have fusion law $\jor{\frac{1}{2}}$ and so $a$ must be a $\jor{\frac{1}{2}}$-axis in $A$. Thus $A_\bt(a)=\{0\}$ and the axet is not skew. 

Suppose $U\cong \widehat{S}(2)^\circ$. This is discussed in Section 5 of \cite{mcinroy2021forbidden}. We have that every non-trivial idempotent in $U$ is a $\jor{\frac{1}{2}}$-axis. As before, we get $A_\bt(a)=\{0\}$ and the axet cannot be skew. 

Therefore $U\cong \TC(\al)$. If $\al\neq-1$, then we get $A\cong \TC(\al,1-\al)$. If $\al= -1$, $A$ is not skew.  \qed
\end{lem}
As we have exhausted all the cases for $P=0$ and $P\neq 0$, this completes the proof of Theorem \ref{Theorem}.

%% file: CorollaryProof.tex
\section{Proof of Corollary \ref{Corollary}}
Let $A=\GenA{p,q}$ where $p$ is a $\mon{\al,\bt}$-axis and $q$ is a $\jor{\al}$-axis. 

Suppose that $q$ is fixed by $\tu{p}$. Both $p,q\in A_{\{0,1,\al\}}(p)$ and by the properties of $\mon{\al,\bt}$, we cannot produce a non-zero element in $A_\bt(p)$. Hence $p$ is a $\jor{\al}$-axis and so $A$ is an axial algebra of $\jor{\al}$-type. 

Now suppose $q$ is not fixed by $p$ and say $q^{\tu{p}}=r$ where $r$ is another axis. By $\tu{p}$ being an automorphism, $r$ is a $\jor{\al}$-axis and $\GenG{p,q}=\{p,q,r\}\cong X'(1+2)$. Therefore $A$ has skew axet of $X'(1+2)$ and is listed in Theorem \ref{Theorem}. 

%% file: Consequences.tex
\section{Larger Skew Axets}
For the majority of this paper, we have set the axet of $(A,X)$ to be $X'(1+2)$. In this section, we will look at larger skew axets and see if our work can simplify their classification.

To remind the reader, let $k\in \N$ and $k>1$. We have the axet of $X'(k+2k)$ to be the axet $X(4k)$ with the extra condition: for $i\in \Z$, $a_{2i}=a_{2(i+k)}$. We do not relabel the vertices when it becomes skew. Therefore,
\[ X(4k):=\{ a_i \; | \; i \in \Z \;\text{ and } \; a_i=a_{i+4k}\}\]
and 
\[ X'(k+2k):=\{ a_i \in X(4k) \; | \; a_{2j}=a_{2(j+k)} \; \forall j \in \Z\}.\]
Looking at the Miyamoto involution of either $X(4k)$ and $X'(k+2k)$, we have that for $i,j\in \Z$, 
\[ a_i^{\tu{j}}= a_{2j-i}.\]
\subsection{The Odd Case}
We will now present a result that M\textsuperscript{c}Inroy suggested.
\begin{prop}\label{odd k}
Let $k$ be odd. For axet $X'(k+2k)$, there exists a subaxet which is isomorphic $X'(1+2)$. 
\proof As $k$ is odd, $a_k$ is in the odd orbit and  $ a_0^{\tu{k}}= a_{2k-0}=a_0$. Therefore $a_0$ is fixed by $\tu{k}$. We have that  $a_k^{\tu{0}}= a_{0-k}=a_{-k}$. Therefore the set $Y=\{a_0,a_k,a_{-k}\}=\GenG{a_0,a_k}$ is a skew axet and is isomorphic to $X'(1+2)$. \qed
\end{prop}
\begin{cor}\label{odd cor}
    Let $k>1$ be odd and $(A,X)$ be a $2$-generated $\mon{\al,\bt}$-axial algebra with skew axet $X'(k+2k)$. Then $\al+\bt=1$. Further, either $A$ has a subalgebra isomorphic to either $\TC(\al,1-\al)$, $Q_2(\frac{1}{3},\frac{2}{3})$ or $Q_2(\frac{1}{3})^\times \oplus \GenG{\id}$.

\proof  By Proposition \ref{odd k}, $X$ must contain a subaxet, $Y=\GenA{a_0,a_k}$, isomorphic to $X'(1+2)$. Let $B=\GenA{a_0,a_k}$ be the $2$-generated axial algebra of Monster type with axet $Y$. As $B$ is generated by two axes in $A$, then $A$ shares the same fusion law in $B$. By Theorem \ref{Theorem}, the rest follows. \qed
\end{cor}

Let $k>1$ and be odd. By Proposition \ref{odd k} and Corollary \ref{odd cor}, we have a possibly strategy to classify all $2$-generated axial algebras of Monster type with skew axet $X'(k+2k)$. Suppose $(A,X)$ has our desired conditions. Notice that $L=\GenA{a_1,a_{-1}}$ is a $2$-generated axial algebra of Monster type. Further, $\tu{0}$ acts as the flip automorphism of $L$ thus $L$ is symmetric. By the classification of $2$-generated symmetric axial algebras of Monster type in \cite{yabe2023classification}, \cite{franchi2022classifying} and \cite{franchi2022quotients} together with the conditions from Corollary \ref{odd cor}, we will have a short list for possible options of $L$. 

We will complete this case in a second paper, going into much more detail. Moreover, we will look at when $k$ is even and progress to completely classifying skew axial algebras of Monster type. 

%% file: Appendix.tex
\appendix
\section{Skew Equations}
We will justify and show the three equations used in Lemma \ref{skew rel} to narrow our search for these skew axial algebras. Although they do not provide much use to understanding how these algebras could be constructed, they do make the proof easier.

Suppose $v$ is an $\mu$-eigenvector of an axis, $x$, where $\mu\neq1$. Then the projection on that axis should be equal to 0; that is, $\lm_x(v)=0$. Coincidentally, nearly all of the eigenvectors in Lemma \ref{eigen a} and \ref{eigen b} satisfy that rule. However we have
\begin{equation*}
 0=\lm_b\left(-\frac{P}{\bt}a+Pb+c\right) = -\frac{P}{\bt}\lmf_1+P+\lmf_2.
\end{equation*}
Whence we get Equation (\ref{proof1}).

\begin{defn}
Let $x$ be a $\mon{\al,\bt}$-axis in $A$, $\lm\in \{1,0, \al, \bt\}$ and $v\in A$. We denote $[v]^x_\lm$ to be the component of $v$ in $ A_\lm(x)$. 
\end{defn}
\begin{lem}
Let $w:=\frac{1}{2}(b-c)$. We have $[a]^a_\bt=0$, $[b]^a_\bt=w$, $[c]^a_\bt=-w$, $[\sg]^a_\bt=0$. Further, $[ab]^a_\bt=\bt w$, $[ac]^a_\bt=-\bt w$, $[bc]^a_\bt=0$, $[a\sg]^a_\bt=0$, $[b\sg]^a_\bt=\dt^fw$, $[c\sg]^a_\bt=-\dt^fw$ and $[\sg^2]^a_\bt=0$.
\end{lem}
\proof
As $a\in A_1(a)$, it has no $\bt$-component in $A_\bt(a)$ and $[a]^a_\bt=0$. As $\sg\in A_{\{1,0,\al\}}(a)$, it has no $\bt$-component in $A_\bt(a)$ and $[\sg]^a_\bt=0$. We can express $b$ in terms of the eigenvectors of $\text{ad}_a$ in Lemma \ref{eigen a}. The reader can check
\[ b= \lm_1 a+ \frac{1}{\al}\left(\ep a+\frac{1}{2}(\al-\bt)(b+c)-\sg\right)+ \frac{1}{\al}\left(\gm a +\frac{1}{2}\bt(b+c)+\sg\right)+\frac{1}{2}(b-c).\]
Thus $[b]_\bt^a=w$. As $c=b^{\tu{a}}$, we get $[c]_\bt^a=-w$.

Let $x, y \in A_{\{0,1,\al\}}(a)$ and notice $x^2, xy\in A_{\{1,0,\al\}}(a)$ and so has no $\bt$-component in $A_\bt(a)$. Therefore $[\sg^2]^a_\bt=[a\sg]^a_\bt=0$. Also
\[ [bc]_\bt^a=P\left([a]_\bt^a+\frac{1}{\bt}[\sg]_\bt^a\right)=0.\]
Note that
\[ [ab]_\bt^a=[\sg]_\bt^a+\bt[a]_\bt^a+\bt[b]_\bt^a=\bt w\]
and 
\[ [b\sg]_\bt^a=(\al-\bt)[\sg]_\bt^a+\bt(\al-\bt)[a]_\bt^a+dt^f[b]_\bt^a=\dt^f w.\]
Applying $\tu{a}$, we get $[ac]_\bt^a$ and $[c\sg]_\bt^a$. \qed

Let $u:= (b -\al)a - \bt b=\sg -(\al-\bt)a$. As $A_\bt(b)=\{0\}$, we have that $u\in A_{\{1,0\}}(b)$. By Lemma \ref{Seress}, the following holds
\[b(au)=(ba)u.\]
Notice
\[ au = a(\sg -(\al-\bt)a)=(\dt -(\al-\bt))a+\frac{1}{2}\bt(\al-\bt)(b+c)+(\al-\bt)\sg\]
and so
\begin{eqnarray*}
[b(au)]_\bt^a &=& (\dt -(\al-\bt))[ab]_\bt^a+\frac{1}{2}\bt(\al-\bt)([b]_\bt^a+[bc]_\bt^a)+(\al-\bt)[b\sg]_\bt^a\\
& =& \left(\bt(\dt -(\al-\bt))+\frac{1}{2}\bt(\al-\bt)+(\al-\bt)\dt^f\right)w
\end{eqnarray*}
We also have 
\begin{eqnarray*}
[(ba)u]_\bt^a&=&[(\sg+\bt a +\bt b)(\sg -(\al-\bt)a)]_\bt^a\\
&=& [\sg^2]_\bt^a -(\al-2\bt)[a\sg]_\bt^a +\bt [b\sg]_\bt^a -\bt(\al-\bt)[a]_\bt^a - \bt(\al-\bt)[ab]_\bt^a\\
&=& (\bt\dt^f -\bt^2(\al-\bt)) w
\end{eqnarray*}
By Lemma \ref{Seress}, we have $0=(ba)u-b(au)$ moreover $0=[(ba)u]_\bt-[b(au)]_\bt$. Looking at the coefficient of $w$, we have
\begin{eqnarray*} 
0&=& (\bt\dt^f-\bt^2(\al-\bt))\\
& -& \left(\bt \dt -\bt(\al-\bt)+\frac{1}{2}\bt(\al-\bt)+(\al-\bt)\dt^f\right)\\
&=&-\bt^2(\al-\bt) -\bt\dt+\frac{1}{2}\bt(\al-\bt)-(\al-2\bt)\dt^f.
\end{eqnarray*}
Rearranging we get Equation (\ref{proof2}).

Let $v:=Pa+\frac{P}{\bt}\sg -\al c=c(b-\al)$. Notice that $v \in A_{\{1,0\}}(b)$. Again by Lemma \ref{Seress}, the following holds
\[b(av)=(ba)v.\]
We have
\begin{eqnarray*}
av &=& Pa +\frac{P}{\bt}\left(\dt a + \frac{1}{2}\bt(\al-\bt)(b+c) +(\al-\bt)\sg\right)\\
& -&\al(\bt a +\bt c +\sg)\\
&=&\left(P +\frac{P}{\bt}\dt -\al\bt\right)a+\left(\frac{1}{2}(\al-\bt)P\right)b\\
&+&\left(\frac{1}{2}(\al-\bt)P-\al\bt\right)c+\left(\frac{P}{\bt}(\al-\bt)-\al\right)\sg.
\end{eqnarray*}
Therefore
\begin{eqnarray*}
[b(av)]_\bt^a &=&\left(P +\frac{P}{\bt}\dt -\al\bt\right)[ab]_\bt^a+\left(\frac{1}{2}(\al-\bt)P\right)[b]_\bt^a\\
&+&\left(\frac{1}{2}(\al-\bt)P-\al\bt\right)[bc]_\bt^a+\left(\frac{P}{\bt}(\al-\bt)-\al\right)[b\sg]_\bt^a.\\
&=&\left(\bt \left(P +\frac{P}{\bt}\dt -\al\bt\right)+\dt^f\left(\frac{P}{\bt}(\al-\bt)-\al\right)\right)w
\end{eqnarray*}
We also have
\begin{eqnarray*}
[(ba)v]_\bt^a&=&\left[\left(\bt a +\bt b +\sg\right)\left(Pa+\frac{P}{\bt}\sg -\al c\right)\right]_\bt^a\\
&=&2P[a\sg]_\bt^a +\frac{P}{\bt}[\sg^2]_\bt^a -\al [c \sg]_\bt^a + \bt P [a]_\bt^a -\al\bt [ac]_\bt^a\\
&+&\bt P [ab]_\bt^a +P[b\sg]_\bt^a -\al\bt [bc]_\bt^a\\
&=&\left(\al \dt^f +\al\bt^2 +\bt^2 P  +P\dt^f\right)w
\end{eqnarray*}
By Lemma \ref{Seress}, $0=[b(av)]^a_\bt-[(ba)v]^a_\bt$ and looking at the coefficient of $w$, we get 
\begin{eqnarray*}
0&=&[b(av)]_\bt-[(ba)v]_\bt\\
&=&\left(\bt P +\dt P -\al\bt^2+\frac{1}{2}(\al-\bt)P+\frac{P}{\bt}(\al-\bt)\dt^f -\al\dt^f\right)\\
&-&\left(\bt^2P +P\dt^f+\al\dt^f +\al\bt^2 \right)\\ 
&=&\left(\frac{P}{\bt}\left[\bt^2 +\bt\dt+\frac{1}{2}\bt(\al-\bt)+(\al-2\bt)\dt^f-\bt^3\right]-2\al(\dt^f+\bt^2)\right).
\end{eqnarray*}
From Equation (\ref{proof2}), we get that
\begin{eqnarray*}
0&=&\frac{P}{\bt}\left[\bt^2 -\bt^2(\al-\bt) +\frac{1}{2}\bt(\al-\bt)-(\al-2\bt)\dt^f\right.\\
&+&\left.\frac{1}{2}\bt(\al-\bt)+(\al-2\bt)\dt^f-\bt^3\right]-2\al(\dt^f+\bt^2)\\
&=&\frac{P}{\bt}\left[\bt^2 -\bt^2(\al-\bt) +\bt(\al-\bt)-\bt^3\right]-2\al(\dt^f+\bt^2)\\
&=&\frac{P}{\bt}\al\bt\left[1-\bt\right]-2\al(\dt^f+\bt^2).
\end{eqnarray*}
Hence we get Equation (\ref{proof3}).

\section*{Acknowledgements}
I would like to thank Professor Sergey Shpectorov for his guidance throughout my PhD studies so far and pushing me to complete this paper. I would also like to thank my family for their continuing support. 